\documentclass[a4paper,11pt]{article}
\usepackage{amsmath,amscd,amssymb}
\usepackage{xypic}
\usepackage[hypertex]{hyperref}
\usepackage{pstricks}
\usepackage{pst-plot}

\def \vs {\vskip}
\def \hs {\hskip}
\def \noi {\noindent}

\def \a {{\alpha}}
\def \b {{\beta}}
\def \ga {{\gamma}}

\def \pd {\varpi}

\def \G {{\mathbb G}}

\def \hole {{\rm Holes}}

\def \Xt {{\widetilde{X}}}

\def \pit {{\widetilde{\pi}}}

\def \Xh {{\widehat{X}}}

\def \pih {{\widehat{\pi}}}

\def \tha #1#2{\noi{\bf#1{\uppercase{\footnotesize{#2}}}}}

\newtheorem{theor}{\tha{T}{heorem}}[section]
\newenvironment{theo}{
  \begin{theor}\hs -0.2 cm {\bf .} ---  }
{  \end{theor}}

\newtheorem{propo}[theor]{\tha{P}{roposition}}
\newenvironment{prop}{
  \begin{propo}\hs -0.2 cm {\bf .} ---  }
{  \end{propo}}

\newtheorem{lemma}[theor]{\tha{L}{emma}}
\newenvironment{lemm}{
  \begin{lemma}\hs -0.2 cm {\bf .} ---  }
{  \end{lemma}}

\newtheorem{fait}[theor]{\tha{F}{act}}
\newenvironment{fact}{
  \begin{fait}\hs -0.2 cm {\bf .} ---  }
{  \end{fait}}

\newtheorem{defini}[theor]{\tha{D}{efinition}}
\newenvironment{defi}{
  \begin{defini}\hs -0.2 cm {\bf .} ---  }
{  \end{defini}}

\newtheorem{corollaire}[theor]{\tha{C}{orollary}}
\newenvironment{coro}{
  \begin{corollaire}\hs -0.2 cm {\bf .} ---  }
{  \end{corollaire}}

\newtheorem{conjecture}[theor]{\tha{C}{onjecture}}
\newenvironment{conj}{
  \begin{conjecture}\hs -0.2 cm {\bf .} ---  }
{  \end{conjecture}}

\newtheorem{exemple}[theor]{\sc{Example}}
\newenvironment{exem}{
  \begin{exemple}\hs -0.2 cm {\bf .} ---  }
{  \end{exemple}}

\newtheorem{remarq}[theor]{\sc{Remark}}
\newenvironment{rema}{
  \begin{remarq}\hs -0.2 cm {\bf .} ---  }
{  \end{remarq}}

\newtheorem{construc}{\sc{Construction}}

\newenvironment{preu}{\noi\sc{Proof}\hs 0.1 cm 
--- \rm }{\hfill$\Box$\vs 0.2 cm}

\def \sca #1#2{\left\langle#1,#2\right\rangle}
\def \scal #1#2{\langle #1,#2 \rangle}

\def \inf {{\preccurlyeq}}
\def \sup {{\succcurlyeq}}

\begin{document}

\topmargin= 1pt

\title{{\bf Gorenstein locus of minuscule Schubert varieties}}
\author{\Large{\sc{Nicolas Perrin}}}
\date{}

\maketitle

\vs -0.5 cm

\begin{small}
\begin{abstract}
In this article, we describe explicitely the Gorenstein locus of all
minuscule Schubert varieties. This proves a special case of a
conjecture of A. Woo and A. Yong \cite{WY} on the Gorenstein locus of
Schubert varieties. 
\end{abstract}
\end{small}

\vs 0.6 cm

\centerline{\large{\textsc{Introduction}}}

\vs 0.3 cm

The description of the singular locus and of the types of
singularities appearing in Schubert varieties is a hard
problem. A first step in this direction was the proof by
V. Lakshmibai and B. Sandhya \cite{LakshSan} of a pattern avoidance
criterion for a Schubert variety in type $A$ to be smooth. There exist
some other results in this direction, for a detailed account see
\cite{LakshBill}. Another
important result was a complete combinatorial description, still in
type $A$, of the irreducible components of the singular locus of a
Schubert variety (this has been realised, almost in the same time, by
L. Manivel \cite{Manivel1} and \cite{Manivel2}, S. Billey and
G. Warrington \cite{BilleyW}, C. Kassel, A. Lascoux and
C. Reutenauer \cite{KasselLR} and A. Cortez \cite{Cortez2}). The
singularity at a generic point of such a component is also given in
\cite{Manivel2} and \cite{Cortez2}. However, as far as I know, this
problem is still open for other types. Another partial result in this
direction is the description of the irreducible components of the
singular locus and of the generic singularity of minuscule and
cominuscule Schubert varieties (see Definition \ref{defmin}) by
M. Brion and P. Polo \cite{Brion-Polo}. 

In the same vein as \cite{LakshSan}, A. Woo and A. Yong gave in
\cite{Alexalex} and \cite{WY} a generalised pattern avoidance
criterion, in type $A$, to decide if a Schubert variety is
Gorenstein. They do not describe the irreducible components of the
Gorenstein locus but give the following conjecture (see Conjecture 6.7
in \cite{WY}):

\begin{conj}
\label{conj}
  Let $X$ be a Schubert variety, a point $x$ in $X$ is in the
  Gorenstein locus of $X$ if and only if the generic point of any
  irreducible component of the singular locus of $X$ containing $x$ is
  is the Gorenstein locus of $X$.
\end{conj}

The interest of this conjecture relies on the fact that, at least in
type $A$, the irreducible components of the singular locus and the
singularity of a generic point of that component are well known. The
conjecture would imply that one only needs to know the information on
the irreducible components of the singular locus to get all the
information on the Gorenstein locus.

In this paper we prove this conjecture for all minuscule Schubert
varieties thanks to a combinatorial description of the
Gorenstein locus of minuscule Schubert varieties. To do this we use
the combinatorial tool introduced in \cite{small} associating to any
minuscule Schubert variety a reduced quiver generalising Young
diagrams. First, we translate the results of M. Brion and P. Polo
\cite{Brion-Polo} in terms of the quiver. We define the holes,
the virtual holes and the essential holes in the quiver (see
Definitions \ref{peak} and \ref{essentiels}) and prove the following:

\begin{theo}
  (\i) A minuscule schubert variety is smooth if and only if its
  associated quiver has no nonvirtual hole.

(\i\i) The irreducible components of the singular locus of a minuscule
Schubert variety are indexed by essential holes.
\end{theo}

Furthermore we explicitely describe in terms of the quiver and the
essential holes these irreducible components and the singularity at a
generic point of a component (for more details see Theorem
\ref{sing}). In particular, with this description it is easy to say if
the singularity at a generic point of an irreducible component of the
singular locus is Gorenstein or not. The essential holes corresponding
to irreducible components having a Gorenstein generic point are called
Gorenstein holes (see also Definition \ref{goren-hole}). We give the
following complete description of the Gorenstein locus:

\begin{theo}
\label{main}
  The generic point of a Schubert subvariety $X(w')$ of a minuscule
  Schubert variety $X(w)$ is in the Gorenstein locus if and only if
  the quiver of $X(w')$ contains all the non Gorenstein holes of the
  quiver of $X(w)$.
\end{theo}

\begin{coro}
\label{maincoro}
  Conjecture \ref{conj} is true for all minuscule Schubert varieties.
\end{coro}

\begin{exem}
Let $\G(4,7)$ be the Grassmannian variety of 4-dimensional
subspaces in a 7-dimensional vector space. Consider the Schubert
variety
$$X(w)=\{V_4\in\G(4,7)\ \dim(V_4\cap W_3)\geq 2\ {\rm and}\ \dim(V_4\cap
W_5)\geq 3\}$$ 
where $W_3$ and $W_5$ are fixed subspaces of dimension 3 and 5
respectively. The minimal length representative $w$ is the permutation
$(2357146)$. Its quiver is the following one (all the arrows are going
down):

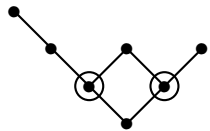
\psset{xunit=1cm}
\psset{yunit=1cm}
\centerline{\begin{pspicture*}(-2,-2)(5,1)
    \psline{-}(0,0)(-0.5,0.5)
    \psline{-}(0,0)(0.5,-0.5)
    \psline{-}(0.5,-0.5)(1,-1)
\psline{-}(1,0)(1.5,-0.5)
\psline{-}(0.5,-0.5)(1,0)
\psline{-}(2,0)(1,-1)
\put(-0.6,0.4){$\bullet$}
\put(-0.1,-0.1){$\bullet$}
\put(0.4,-0.6){$\bullet$}
\put(0.9,-0.1){$\bullet$}
\put(1.4,-0.6){$\bullet$}
\put(0.9,-1.1){$\bullet$}
\put(1.9,-0.1){$\bullet$}
\pscircle(0.5,-0.5){0.2}
\pscircle(1.5,-0.5){0.2}
\end{pspicture*}}

We have circled the two holes on this quiver. The left hole is not a
Gorenstein hole (this can be easily seen because the two peaks above
this hole do not have the same height, see Definition \ref{peak}) but
the right one is Gorenstein (the two peaks have the same height). Let
$X(w')$ be an irreducible component of the singular locus of
$X(w)$. The possible quivers of such a variety $X(w')$ are the
following (for each hole we remove all the vertices above that hole):

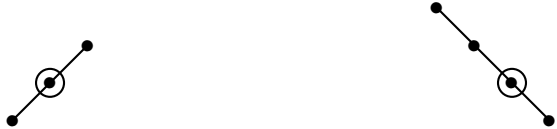
\psset{xunit=1cm}
\psset{yunit=1cm}
\centerline{\begin{pspicture*}(-2,-2)(5,1)
\psline{-}(2,0)(1,-1)
\put(1.4,-0.6){$\bullet$}
\put(0.9,-1.1){$\bullet$}
\put(1.9,-0.1){$\bullet$}
\pscircle(1.5,-0.5){0.2}
\end{pspicture*}
\begin{pspicture*}(-2,-2)(5,1)
    \psline{-}(0,0)(-0.5,0.5)
    \psline{-}(0,0)(0.5,-0.5)
    \psline{-}(0.5,-0.5)(1,-1)
\put(-0.6,0.4){$\bullet$}
\put(-0.1,-0.1){$\bullet$}
\put(0.4,-0.6){$\bullet$}
\put(0.9,-1.1){$\bullet$}
\pscircle(0.5,-0.5){0.2}
\end{pspicture*}}

\noi
These Schubert varieties correspond to the permutations: $(1237456)$
and $(2341567)$. Let
$X(w')$ be a Schubert subvariety in $X(w)$ whose generic point is not
in the Gorenstein locus. Then $X(w')$ has to be contained in
$X(1237456)$.
\end{exem}

\vs 0.2 cm

\noi
{\bf Acknowledgements:} I thank Frank Sottile and Jim Carrel for their
invitation to the BIRS workshop \emph{Comtemporary Schubert calculus}
during which the major part of this work has been done. 

\section{Minuscule Schubert varieties}

Let us fix some notations and recall the definitions of minuscule
homogeneous spaces and minuscule Schubert varieties. A basic
reference is \cite{GofG/P3}.

In this paper $G$ will be a semi-simple algebraic group, we fix
$B$ a Borel subgroup and $T$ a maximal torus in $B$. We denote by $R$
the set of roots, by $R^+$ and $R^-$ the set of positive and negative
roots. We denote by $S$ the set of simple roots. We will denote by $W$
the Weyl group of $G$.

We also fix $P$ a parabolic subgroup containing $B$. We denote by
$W_P$ the Weyl group of $P$ and by $W^P$ the set of minimal length
representatives in $W$ of the coset $W/W_P$. Recall that the Schubert
varieties in $G/P$ (that is to say the $B$-orbit closures in $G/P$)
are parametrised by $W^P$.

\begin{defi}
A fundamental weight $\pd$ is said to be minuscule if, for all
positive roots $\a\in R^+$, we have $\sca{\a^\vee}{\pd}\leq1$.
\end{defi}
 
With the notation of N. Bourbaki \cite{bourb}, the minuscule weights
are:

\begin{center}
\begin{tabular}{|c|c|c|}
\hline
Type&minuscule\\
\hline
$A_n$&$\pd_1\cdots\pd_n$\\
\hline
$B_n$&$\pd_n$\\
\hline
$C_n$&$\pd_1$\\
\hline
$D_n$&$\pd_1$, $\pd_{n-1}$ and $\pd_n$\\
\hline
$E_6$&$\pd_1$ and $\pd_6$\\
\hline
$E_7$&$\pd_7$\\
\hline
$E_8$&none\\
\hline
$F_4$&none\\
\hline
$G_2$&none\\
\hline
\end{tabular}
\end{center}

\begin{defi}
  \label{defmin}
  Let $\pd$ be a  minuscule weight and let $P_\pd$ be the associated
  parabolic subgroup. The homogeneous space $G/P_\pd$ is then said to
  be minuscule. The Schubert varieties of a minuscule homogeneous
  space are called minuscule Schubert varieties.
\end{defi}

\begin{rema}
  It is a classical fact that
to study minuscule homogeneous spaces and their Schubert varieties, it
is sufficient to restrict ourselves to simply-laced groups.
\end{rema}

In the rest of the paper, the group $G$ will be simply-laced, the
subgroup $P$ will be a
maximal parabolic subgroup associated to a minuscule fundamental
weight $\pd$. The minuscule homogeneous space $G/P$ will be denoted by
$X$ and the Schubert variety associated to $w\in W^P$ will be denoted
by $X(w)$ with the convention that the dimension of $X(w)$ is the
length of $w$.

\section{Miniscule quivers}

In \cite{small}, we associated to any minuscule Schubert variety
$X(w)$ a unique quiver $Q_w$. The definition a priori depends on the
choice of a reduced expression but does not depend on the commuting
relations. In the minuscule setting this implies that the following
definitons do not depend on the choosen reduced expression. Fix a
reduced expression
$w=s_{\b_1}\cdots s_{\b_r}$
of $w$ (recall that $w$ is in $W^P$ the set of minimal length
representatives of $W/W_P$) where for all $i\in[1,r]$, we have
$\b_i\in S$.

\begin{defi}
\label{quiver}
(\i) The successor $s(i)$ and the predecessor $p(i)$ of an element
$i\in[1,r]$ are the elements
$\displaystyle{s(i)=\min\{j\in[1,r]\ /\ j>i\ \textrm{\textit{and}}\
  \b_j=\b_i\}}$ and $\displaystyle{p(i)=\max\{j\in[1,r]\ /\
  j<i\ \textrm{\textit{and}}\ \b_j=\b_i\}}.$

(\i\i) Denote by $Q_w$ the quiver whose set of vertices is the set
  $[1,r]$ and whose arrows are given in the 
  following way: there is an arrow from $i$ to $j$ if and only if 
  $\scal{\b_j^\vee}{\b_i}\neq0$ and $i<j<s(i)$ (or only $i<j$ if $s(i)$
  does not exist).
\end{defi}

\begin{rema}  
(\i) This quiver comes with a coloration of its vertices by simple roots
via the map $\b:[1,r]\to S$ such that $\b(i)=\b_i$.

(\i\i) There is a natural order on the quiver $Q_w$ given by $i\inf j$
if there is an oriented path from $j$ to $i$. Caution that this order
is the reversed order of the one defined in \cite{small}.

(\i\i\i) Note that if we denote by $Q_\pd$ the quiver obtained from
the longuest element in $W^P$, then the quiver $Q_w$ is a subquiver of
$Q_\pd$. The quivers of Schubert subvarieties are exactely the
order ideals in the quiver $Q_\pd$. We will call such a quiver reduced
(meaning that it corresponds to a reduced expression of an element in
$W^P$, see \cite{small} for more details on the shape of reduced
quivers). 
\end{rema}

Recall also that we defined in \cite{small} some combinatorial objects
associated to the quiver $Q_w$.

\begin{defi}
\label{peak}
(\i) We call peak any vertex of $Q_w$ maximal for the partial order
  $\preccurlyeq$. We denote by ${\rm Peaks}(Q_w)$ the set of peaks of $Q_w$.

(\i\i) We call hole of the quiver $Q_w$ any vertex $i$ of $Q_\pd$
satisfying one of the following properties
\begin{itemize}
\item the vertex $i$ is in $Q_w$ but $p(i)\not\in Q_w$
 and there are exactly two
vertices $j_1\sup i$ and $j_2\sup i$ in $Q_w$ with
$\sca{\b_i^\vee}{\b_{j_k}}\neq0$ for $k=1,2$.
\item the vertex $i$ is not in $Q_w$, $s(i)$ does not exist in $Q_\pd$
  and there exist $j\in Q_w$ with $\sca{\b_i^\vee}{\b_j}\neq0$. 
\end{itemize}
Because the vertex of the second type of holes is not a vertex in
$Q_w$ we call such a hole a virtual hole of $Q_w$.
We denote by $\hole(Q_w)$ the set of holes of $Q_w$.

(\i\i\i) The height $h(i)$ of a vertex $i$ is the largest positive
integer $n$ such that
there exists a sequence $(i_k)_{k\in[1,n]}$ of vertices with $i_1=1$,
$i_n=r$ and such that there is an arrow from $i_k$ to $i_{k+1}$ for
all $k\in[1,n-1]$.
\end{defi}

Many geometric properties of the Schubert variety $X(w)$ can be read
on its quiver. In particular we proved in \cite[Corollary 4.12]{small}:

\begin{prop}
\label{stab}
  A Schubert subvariety $X(w')$ in $X(w)$ is stable under ${\rm
    Stab}(X(w))$ if and only if $\b({\rm Holes}(Q_{w'}))\subset\b({\rm
    Holes}(Q_w))$.
\end{prop}

An easy consequence of this fact and the result by M. Brion and
P. Polo that the smooth locus of $X(w)$ is the dense ${\rm
  Stab}(X(w))$-orbit is the following:

\begin{prop}
  A Schubert variety $X(w)$ is smooth if and only if all the holes of
  its quiver $Q_w$ are virtual.
\end{prop}

We will be more precise in Theorem \ref{sing} and we will describe the
irreducible components of the singular locus and the generic
singularity of this component in terms of the quiver. The
Gorensteiness of the variety is also easy to detect on the quiver as
we proved in \cite[Corollary 4.19]{small}:

\begin{prop}
  A Schubert variety $X(w)$ is Gorenstein if and only if all the peaks
  of its quiver $Q_w$ have the same height.
\end{prop}

\section{Generic singularities of minuscule Schubert varieties}

In this section, we go one step further in the direction of reading on
the quiver $Q_w$ the geometric properties of $X(w)$. We will translate
the results of M. Brion and P. Polo \cite{Brion-Polo} on the
irreducible components of the singular locus of $X(w)$ and the
singularity at a generic point of such a component in terms of the
quiver $Q_w$. We will need the following notations:

\begin{defi}
\label{essentiels} 
(\i) Let $i$ be a vertex of $Q_w$, we define the subquiver $Q_{w}^i$
of $Q_w$ as the full subquiver containing the following set of
vertices $\{j\in Q_w\ /\ j\succcurlyeq i\}.$ We denote by $Q_{w,i}$
the full subquiver of $Q_w$ containing the vertices of $Q_w\setminus
Q_w^i$. We denote by $w^i$ (resp. $w_i$) the elements in $W^P$
associated to the quivers $Q_w^i$ (resp. $Q_{w,i}$).

(\i\i) A hole $i$ of the quiver $Q_w$ is said to be essential if it is
not virtual and if there is no hole in the subquiver $Q_{w}^i$. 

(\i\i\i) Following M. Brion and P. Polo, denote by $J$ the set
$\beta({\rm Holes}(Q_w))^c$.
\end{defi}

We then prove the following:

\begin{theo}
\label{sing}
(\i) The set of irreducible components of the singular locus of $X(w)$
is in one to one correspondence with the set of essential holes of the
quiver $Q_w$.
In particular, if $i$ is an essential hole of $Q_w$, the corresponding
irreducible component is the Schubert subvariety $X(w_i)$ of $X(w)$
whose quiver is $Q_{w,i}$.

(\i\i\i) Furthermore, the singularity of $X(w)$ at a generic point of
$X(w_i)$ is the same singularity as the one of the $B$-fixed point in
the Schubert variety $X(w^i)$ whose quiver is $Q_w^i$.
\end{theo}

\begin{rema}
The singularity of the
  $B$-fixed point in $X(w^i)$ is described in
  \cite{Brion-Polo}.
\end{rema}

  \begin{preu}
This result is a reformulation of the main results of M. Brion and
P. Polo \cite{Brion-Polo}.
Proposition \ref{stab} shows that the essential holes are in one to
one correspondence with maximal Schubert subvarieties in $X(w)$ stable
under ${\rm Stab}(X(w))$ and that if $i$ is an essential hole, then
the corresponding Schubert subvariety $X(w_i)$ is associated to the
quiver $Q_{w,i}$. According to \cite{Brion-Polo}, these are the
irreducible components of the singular locus.

To describe the singularity of $X(w_i)$, M. Brion and P. Polo define
two subsets $I$ and $I'$ of the set of simple roots as follows:
\begin{itemize}
\item the set $I$ is the union of the connected components of $J\cap
  w_i(R_P)$ adjacent to $\b(i)$ 
\item the set $I'$ is the union $I\cup\{\b(i)\}$.
\end{itemize}
We describe these sets thanks to the quiver.

\begin{prop}
\label{I'}
  The set $I'$ is $\b(Q_w^i)$.
\end{prop}

\begin{preu}
The elements in $J\cap w_i(R_P)$ are the simple roots $\ga\in J$
such that ${w_i}^{-1}(\ga)\in R_P$. Thanks to Lemma \ref{im-rac}, these
elements are the simple roots in $J$ neither in $\b({\rm
  Holes}(Q_{w,i}))$ nor in $\b({\rm Peaks}(Q_{w,i}))$. 

An easy (but fastidious for types $E_6$ and $E_7$) look on the quivers
shows that $I'=\b(Q_w^i)$. A uniform proof of this statement is
possible but needs an involved case analysis on the quivers. 
    \end{preu}

\begin{lemm}
\label{im-rac}
Let $\beta$ be a simple root, then we have
\begin{itemize}
 \item[1.] $w^{-1}(\b)\in R^-\setminus R^-_P$ if $\b\in \b({\rm Peaks}(Q_w))$,
 \item[2.] $w^{-1}(\b)\in R^+\setminus R^+_P$ if $\b\in \b({\rm
     Holes}(Q_w))=J^c$ or
 \item[3.] $w^{-1}(\b)\in R^+_P$ otherwise.
\end{itemize}
\end{lemm}

\begin{preu}
Let $w=s_{\b_1}\cdot s_{\b_r}$ be a reduced expression for $w$, we
want to compute $w^{-1}(\b)=s_{\b_r}\cdots s_{\b_1}(\b)$. We proceed
by induction and deal with the three cases at the same time.

1. Take first $\beta\in \b({\rm Peaks}(Q_w))$, we may assume that $\b_1=\b$ and
$w^{-1}(\b)=s_{\b_r}\cdots s_{\b_2}(-\b)$. Let $i\in{\rm Peaks}(Q_w)$
such that $\beta(i)=\b$, the quiver obtained by removing $i$ has
$s(i)$ for hole (possibly virtual). We may apply induction and the
result in case 2.

2.a. Let $\b\in J^c$. Assume first that there is no $k\in Q_w$ with
$\b(k)=\b$. Then there exist an $i\in Q_w$ such that
$\sca{\b^\vee}{\b_i}\neq0$. Let us prove that such a vertex $i$ is
unique. Indeed, the support of $w$ is contained in a subdiagram $D$
of the Dynkin diagram not containing $\b$. The diagram $D$ contains
the simple root $\a$ corresponding to $P$ (except if $X(w)$ is a point
in which case $w={\rm Id}$ and the lemma is easy). The quiver $Q_w$ is
in particular contained in the quiver of the minuscule homogeneous
variety associated to $\a\in D$. It is easy to check on these
quivers (see in \cite{small} for the shape of these quivers) that
there is a unique such vertex $i$. 

Now consider the quivers $Q_w^i$ and $Q_{w,i}$. Recall that we denote
by $w^i$ and $w_i$ the associated elements in $W$. We have
$w=w^iw_i$. We compute ${w^i}^{-1}(\b)$ and because all simple roots
$\b(x)$ for $x\in Q_w^i$ with $x\neq i$ are orthogonal to $\b$ we have
${w^i}^{-1}(\b)=s_{\b_i}(\b)=\b+\b_i$. We then have
${w}^{-1}(\b)=w_i^{-1}(\b+\b_i)$. Because $i$ was the only vertex such
that $\sca{\b^\vee}{\b_i}\neq0$, we have $w_i^{-1}(\b)=\b\in R_P^+$
and by induction (note that $i$ is now a hole of $Q_{w,i}$) we have
$w_i^{-1}(\b_i)\in R^+\setminus R_P^+$ and we have the result. 

2.b. Now assume that there exist $k\in {\rm Holes}(Q_w)$ with
$\b(k)=\b$ and let $i$ a vertex maximal for the property
$\sca{\b^\vee}{\b_i}\neq0$. Remark that we have $k<i$. Consider one
more time the quivers $Q_w^i$ and $Q_{w,i}$ and the elements $w^i$ and
$w_i$. We have $w^{-1}(\b)=w_i^{-1}(\b_i+\b)$. But as before we have
by induction $w_i^{-1}(\b_i)\in R^+\setminus R_P^+$ so that we can
conclude by induction as soon as $k$ is not a peak of $Q_{w,i}$. But
because $k$ is an hole, there exist a vertex $j\in Q_w$ with $j\neq i$
and such that there is an arrow $j\to k$ in $Q_w$. Because $i$ was
taken maximal $j$ is a vertex of $Q_{w,i}$ and $k$ is not a peak of
this quiver. 

3. If $\b$ is not in the support of $w$ but is not in $\b({\rm
  Holes})$ then $w^{-1}(\b)=\b\in R^+_P$. 

Let $\b$ in $\beta(Q_w)$ but not in $\b({\rm Holess}(Q_w))$ or
$\b({\rm Peaks}(Q_w))$ and let $k$ the highest 
vertex such that $\b(k)=\b$. There exists a unique vertex $i\in Q_w$
such that $i\succ k$ and $\sca{\b^\vee}{\b(i)}\neq0$. We have
$w^{-1}(\b)=w_i^{-1}(\b_i+\b)$ and the vertex $k$ is a peak of
$Q_{w,i}$ so that $w_i=s_{\b(k)}w_k=s_{\b}w_k$ and
$w^{-1}(\b)=w_k^{-1}(\b_i)$. Now it is easy to see that either $s(i)$
does not exists and in this case it is not a virtual hole or it exists
but is neither a peak nor a hole of $Q_{w,k}$. We conclude by
induction on the third case.
\end{preu}

The Theorem is now a corollary of the description of the singularities
thanks to $I$ and $I'$ done by M. Brion and P. Polo.
\end{preu}

\begin{rema}
In their article M. Brion and P. Polo also deal with the cominucule
Schubert varieties. We believe that, in that case, Theorem \ref{main}
should hold true as well as Corollary \ref{maincoro}.
\end{rema}

It is now easy to decide which generic singularity is Gorenstein:

\begin{coro}
Let $i$ be an essential hole of the quiver $Q_w$. The generic point of
the irreducible component $X(w_i)$ of the singular locus is Gorenstein
if and only if all the peaks of $Q_w^i$ are of the same height. 
\end{coro}

We describe the Schubert subvarieties $X(w')$ in
$X(w)$ that are expected to be Gorenstein at their generic point by
the conjecture of A. Woo and A. Yong. Let us give the following

\begin{defi}
\label{goren-hole}
(\i) An essential hole is said to be Gorenstein if the generic point
of the associated irreducible component of the singular locus is in
the Gorenstein locus.

(\i\i) A Schubert subvariety $X(w')$ in $X(w)$ is said to have the property
  \emph{(WY)} if the generic point of any irreductible component of the
  singular locus of $X(w)$ containing $X(w')$ is in the Gorenstein
  locus of $X(w)$.
\end{defi}

We have the following:

\begin{prop}
\label{sensfacile}
Let $X(w')$ be a Schubert subvariety of the Schubert variety
$X(w)$. If the generic point of $X(w')$ is Gorentein in $X(w)$, then
$X(w')$ has the property \emph{(WY)}.
\end{prop}

\begin{preu}
  Let $X(v)$ be an irreducible component of the singular locus of
  $X(w)$ containing $X(w')$. Because the property of beeing non
  Gorenstein is stable under closure, this implies that the generic
  point of $X(v)$ is Gorenstein in $X(w)$. 
\end{preu}

Remark that, because all the irreducible components of the singular
locus of $X(w)$ are stable under ${\rm Stab}(X(w))$, the property
(WY) need only to be checked on ${\rm Stab}(X(w))$-stable Schubert
subvarieties.

\begin{prop}
\label{lieuWY}
  (\i) The Schubert subvarieties $X(w')$ in $X(w)$ stable under ${\rm
    Stab}(X(w))$ are exactely those such that the associated quiver
  $Q_{w'}$ satisfies
$$Q_ {w'}=\bigcap_{i\in{\rm Holes}(Q_w)}Q_{w,s^{k_i}(i)}$$
where the $(k_i)_{i\in{\rm Holes}(Q_w)}$ are integers greater or equal
to $-1$ (if $k_i=-1$, the quiver $Q_{w,s^{k_i}(i)}$ is $Q_w$ by
definition).

(\i\i) A ${\rm Stab}(X(w))$-stable Schubert subvariety $X(w')$ of
$X(w)$ has the property \emph{(WY)} if and only if the only essential
holes in the difference $Q_w\setminus Q_{w'}$ are
Gorenstein. Equivalentely, writing 
$$Q_ {w'}=\bigcap_{i\in{\rm Holes}(Q_w)}Q_{w,s^{k_i}(i)},$$
if and only if the only holes in of the quivers
$(Q_{w}^{s^{k_i}(i)})_{i\in{\rm Holes}(Q_w)}$ are Gorenstein
holes. Another equivalent formulation is that $Q_{w'}$ contains all the
non Gorenstein essential holes of $Q_w$.
\end{prop}

\begin{preu}
  (\i) Consider the subquiver $Q_{w'}$ in $Q_w$ and for each hole $i$
  of $Q_w$ define the integer $k_i=\min\{k\geq0\ /\
  s^k(i)\in Q_{w'}\}-1$. Because of the fact (see for example
  \cite{GofG/P3}) that the strong and weak Bruhat orders coincide for
  minuscule Schubert varieties, the quiver $Q_{w'}$ has to be
  contained in the intersection
$$Q'=\bigcap_{i\in{\rm Holes}(Q_w)}Q_{w,s^{k_i}(i)}.$$
We therefore need to remove some vertices to $Q'$ to get $Q_{w'}$. But
removing a vertex $j$ of the quiver $Q'$ (it has to be a peak of
$Q'$) creates a hole in $s(j)$ (or a virtual hole in $j$ if $s(j)$
does not exist). Because $X(w')$ is ${\rm Stab}(X(w))$-stable, the
last removed vertex $j$ is such that $\b(j)\in\b({\rm
  Holes}(Q_w))$. This implies that no more vertex can be removed from
$Q'$ to get $Q_{w'}$ and in particular $Q_{w'}=Q'$.

(\i\i) The Schubert subvariety has the property (WY) if and only if
all the irreducible components $X(w_i)$ of the singular locus of
$X(w)$ containing $X(w')$ are such that $i$ is a Gorenstein hole. But
$X(w')$ is contained in $X(w_i)$ if and only if $Q_{w'}$ is contained
in $Q_{w,i}$. This is equivalent to the fact that $Q_w^i$ is contained
in $Q_w\setminus Q_{w'}$ and the proof follows. 
\end{preu}

\section{Relative canonical model and Gorenstein locus}

In this section, we recall the explicit construction given in
\cite{small} of the relative canonical model of $X(w)$. Recall that we
described in \cite{small} the Bott-Samelson resolution $\pi:\Xt(w)\to
X(w)$ as a configuration variety {\`a} la Magyar \cite{Magyar-conf1}:
$$\Xt(w)\subset\prod_{i\in Q_w}G/P_{\b_i}$$
where $P_{\b_i}$ is the maximal parabolic associated to the simple
root $\b_i$. The map $\pi:\Xt(w)\to X(w)$ is given by the projection 
$\displaystyle{\prod_{i\in Q_w}G/P_{\b_i}\to G/P_{\b_{m(w)}}}$ where
$m(w)$ is the smallest element in $Q_w$.

\vs 0.2 cm

We define a partition on the peaks of the quiver $Q_w$ and a partition
of the quiver itself:

\begin{defi}
  \label{lacanonique}
(\i) Define a partition $(A_i)_{i\in[1,n]}$ of ${\rm Peaks}(Q_w)$ by
induction: $A_1$ is the 
  set of peaks with minimal height and $A_{i+1}$ is the set of peaks in
  ${\rm Peaks}(Q_w)\setminus\bigcup_{k=1}^iA_k$ with minimal height
  (the integer $n$ is the number of different values the height
  function takes on the set ${\rm Peaks}(Q_w)$).

(\i\i) Define a partition $(Q_w(i))_{i\in[1,n]}$ of $Q_w$ by induction:
$$Q_w(i)=\{x\in Q_w\ /\ \exists j\in A_i\ : x\preccurlyeq j\ {\rm
  and}\ x\not\preccurlyeq k\ \forall k\in\cup_{j>i}A_j\}.$$
\end{defi}

We proved in \cite{small} that these quivers $Q_w(i)$ are quivers of
minuscule Schubert varieties and in particular have a minimal element
$m_w(i)$. We defined the variety $\Xh(w)$ as the image of the
Bott-Samelson resolution $\Xt(w)$ (seen as a configuration variety)
in the product $\prod_{i=1}^nG/P_{\b_{m_w(i)}}$. 

Because $m_w(n)=m(w)$
we have a map $\pih:\Xh(w)\to X(w)$ and a factorisation
$$\xymatrix{\Xt(w)\ar[dr]_\pi\ar[r]^\pit&\Xh(w)\ar[d]^\pih\\
&X(w).}$$

We proved the following result in \cite{small}:

\begin{theo}
  (\i) The variety $\Xh(w)$ together with the map $\pih$ realise
  $\Xh(w)$ as the relative canonical model of $X(w)$.

(\i\i) The variety $\Xh(w)$ is a tower of locally trivial fibrations
with fibers the Schubert varieties associated to the quivers
$Q_w(i)$. In particular $\Xh(w)$ is Gorenstein.
\end{theo}

We will use this resolution to prove our main result. Indeed, we will
prove that the generic fibre of the map $\pih:\Xh(w)\to X(w)$ above a
(WY) Schubert subvariety $X(w')$  is a point. In other words, the map
$\pih$ is an isomorphism on an open subset of $X(w')$. As a
consequence, the generic point of $X(w')$ will be in the Gorenstein
locus.

Let us recall some facts on $\Xt(w)$ and $\Xh(w)$ (see \cite{small}):

\begin{fact}
\label{imageBS}
  (\i) To each vertex $i$ of $Q_w$ one can associated a divisor $D_i$
  on $\Xt(w)$ and all these divisors intersect transversally. 

(\i\i) For $K$ a subset of the vertices of $Q_w$, we denote by $Z_K$
the transverse intersection of the $D_i$ for $i\in K$.

(\i\i\i) The image of the closed subset $Z_K$ by the map $\pi$ is the
Schubert variety $X(w_K)$ whose quiver $Q_{w_K}$ is the biggest
reduced subquiver of $Q_w$ not containing the vertices in $K$.
\end{fact}

The quiver $Q_w(i)$ defines a element $w(i)$ in $W$ and the fact that
these quivers realise a partition of $Q_w$ implies that we have an
expression $w=w(1)\cdots w(n)$ with $l(w)=\sum l(w(i))$. We prove the
following generalisation of this fact:

\begin{prop}
\label{image}
Let $K$ be a subset of the vertices of $Q_w$. The image of the closed
subset $Z_K$ by the map $\pit$ is a tower of locally trivial
fibrations with fibers the Schubert varieties $X(w_K(i))$ whose quiver
$Q_{w_K(i)}$ is the biggest reduced subquiver of $Q_{w(i)}$ not
containing the vertices of $K\cap Q_{w(i)}$. 

This variety is the image by $\pit$ of $Z_{\cup_{i=1}^nQ_K(i)}$.
\end{prop}

\begin{preu}
As we explained in \cite[Proposition 5.9]{small}, the Bott-Samelson
resolution is the quotient of the product $\prod_{i\in Q_w} R_i$ where
the $R_i$ are certain minimal parabolic subgroups by a product of
Borel subgroups $\prod_{i=1}^r B_i$. The variety $\Xh(w)$ is the
quotient of a product $\prod_{i=1}^nN_i$ of parabolic subgroups such
that the multiplication in $G$ maps $\prod_{k\in Q_{w(i)}}R_k$ to
$N_i$ by a product $\prod_{i=1}^nM_i$ of parabolic subgroups. The map
$\pit$ is induced by the product from $\prod_{i\in Q_w} R_i$ to
$\prod_{i=1}^nN_i$. In particular, this means that for $i\in[1,n]$
fixed, the map $\prod_{k\in Q_{w(i)}}\to N_i$ induces the map from the
Bott-Samelson resolution $\Xt(w(i))$ to $X(w(i))$. We may now apply
part \textit{(\i\i\i)} of the preceding fact because the quiver
$Q_w(i)$ is minuscule.
\end{preu}

We now remark that the quivers $Q_{w'}$ associated to Schubert
subvarieties $X(w')$ in the Schubert variety $X(w)$ having the
property (WY) have a nice behaviour with repect to the partition
$(Q_w(i))_{i\in [1,n]}$ of $Q_w$.

\begin{prop}
  Let $X(w')$ be a ${\rm Stab}(X(w))$-stable Schubert subvariety of
  $X(w)$ having the property \emph{(WY)}. Let us denote by
  $(C_j)_{j\in[1,k]}$ the connected components of the subquiver
  $Q_w\setminus Q_{w'}$ of $Q_w$. Then for each $j$, there exist an
  unique $i_j\in[1,n]$ such that $C_j\subset Q_w(i_j)$.
\end{prop}

\begin{preu}
Recall from Proposition \ref{lieuWY} that, denoting by ${\rm
  GorHol}(Q_w)$ the set of Gorenstein holes in $Q_w$, we may write 
$$Q_w\setminus Q_ {w'}=\bigcup_{i\in{\rm GorHol}(Q_w)}Q_{w}^{s^{k_i}(i)}$$
with $k_i$ an integer greater or equal to $-1$ and with the additional
condition that $Q_{w}^{s^{k_i}(i)}$ contains only Gorenstein holes. 
Because the quivers $Q_{w}^{s^{k_i}(i)}$ are connected, any connected
component of $Q_w\setminus Q_ {w'}$ is an union of such quivers. But
we have the following:

\begin{lemm}
  Let $i\in{\rm Holes}(Q_w)$ and assume that $Q_w^{s^k(i)}$ meets at
  least two subquivers of the partition $(Q_w(i))_{i\in [1,n]}$, then
  $Q_w^{s^k(i)}$ contains a non Gorenstein hole.
\end{lemm}

\begin{preu}
  The quiver $Q_w^{s^k(i)}$ meets two subquivers of the partition
  $(Q_w(i))_{i\in [1,n]}$, in particular it contains two peaks of
  $Q_w$ of different heights. By connexity of $Q_w^{s^k(i)}$, we may
  assume that these two peaks are adjacent. In particular there is a
  hole between these two peaks and this hole is not Gorenstein and is
  contained in $Q_w^{s^k(i)}$.
\end{preu}

The proposition follows.
\end{preu}

We describe the inverse image by $\pih$ of a ${\rm
  Stab}(X(w))$-stable Schubert subvariety of $X(w)$ having the
property (WY). To do this, first remark that the map $\pi$ is
$B$-equivariant and that the inverse image $\pi^{-1}(X(w'))$ has to be
a union of closed subsets $Z_K$ for some subsets $K$ of $Q_w$. Let
$Z_K\subset \pi^{-1}(X(w'))$ be such that $\pi:Z_K\to X(w')$ is
dominant. We will denote by $Q_w^{w'}(i)$ the intersection $Q_{w'}\cap
Q_w(i)$ and by $w'(i)$ the associated element in $W$.

\begin{prop}
  The image of $Z_K$ in $\Xh(w)$ by $\pit$ is the same as the image of
  $Z_{Q_w\setminus Q_{w'}}$.
\end{prop}

\begin{preu}
  Thanks to Proposition \ref{image} we only need to compute the
  quivers $Q_{w_K(i)}$. Consider the decomposition into connected
  components $Q_w\setminus Q_{w'}=\cup_{j=1}^k C_j$. We may decompose
  $K$ accordingly as $K=\cup_{j=1}^kK_j$ where $K_j=K\cap C_j$. But
  because each connected component of $Q_w\setminus Q_{w'}$ is
  contained in one of the quivers $(Q_w(i))_{i\in[1,n]}$ this implies
  that $Q_{w_K(i)}$ is exactely $Q_{w_K}\cap Q_w(i)$ where $Q_{w_K}$
  is the biggest reduced quiver in $Q_w$ $Q_w$ not containing the
  vertices in $K$ (see Fact \ref{imageBS}). We get $Q_{w_K}=Q_{w'}$
  (because $Z_K$ is sent onto $X(w')$) and the result follows.
\end{preu}

\begin{theo}
  Let $X(w')$ be a Schubert subvariety in $X(w)$. Then $X(w')$ has the
  property \emph{(WY)} if and only if its generic point is in the Gorenstein
  locus of $X(w)$.
\end{theo}

\begin{preu}
  We have already seen in Proposition \ref{sensfacile} that if the
  generic point of $X(w')$ is in the Gorenstein locus of $X(w)$ then
  $X(w')$ has the property (WY). 

Conversely let $X(w')$ be a Schubert subvariety having the property
(WY). The previous proposition implies that its inverse image
$\pih^{-1}(X(w'))$ is the variety $\pit(Z_{Q_w\setminus Q_{w'}})$. But
this implies that the map $\pih:\pit(Z_{Q_w\setminus
  Q_{w'}})=\pih^{-1}(X(w'))\to 
X(w')$ is birational (because the varieties have the same dimension
given by the number of vertices in the quiver). In particular, the map
$\pih$ is an isomorphism on an open subset of $X(w)$ meeting $X(w')$
non trivially. Therefore, because $\Xh(w)$ is Gorenstein, it is the
case of the generic point in $X(w')$ as a point in $X(w)$. 
\end{preu}

\begin{small}

\vs 0.2 cm

\noi
{\textsc{Universit{\'e} Pierre et Marie Curie - Paris 6}}

\vs -0.1 cm

\noi
{\textsc{UMR 7586 --- Institut de Math{\'e}matiques de Jussieu}}

\vs -0.1 cm

\noi
{\textsc{175 rue du Chevaleret}}

\vs -0.1 cm

\noi
{\textsc{75013 Paris,}} \hs 0.2 cm{\textsc{France.}}

\vs -0.1 cm

\noi
{email : \texttt{nperrin@math.jussieu.fr}}

\end{small}

\end{document}